\documentclass[singlespacing]{elsart}

\usepackage{amsfonts,amsmath,amssymb}
\usepackage{slashed,setspace}

\usepackage{enumerate}
\usepackage[active]{srcltx}
\def\leftnote#1{\vadjust{\setbox1=\vtop{\hsize 20mm\parindent=0pt\bf\baselineskip=9pt\rightskip=4mm plus 4mm#1}\hbox{\kern-25mm\smash{\box1}}}}
\newtheorem{ithm}{Theorem}

\newtheorem{theorem}{Theorem}[section]

\newtheorem{lemma}[theorem]{Lemma}
\newtheorem{remark}[theorem]{Remark}
\newtheorem{corollary}[theorem]{Corollary}
\newtheorem{proposition}[theorem]{Proposition}

\newenvironment{proof}[1][Proof]{\par\addvspace{2mm}\noindent\textbf{#1.} }{\ \rule{0.5em}{0.5em}\par\vspace{4mm}}

\newcommand{\wb}{\overline}

\newcommand{\Z}{\mathbb{Z}}
\newcommand{\FTI}{\overline{\cal F}}
\newcommand{\FT}{\mathcal{F}}
\newcommand{\F}{\mathbb{F}}

\newenvironment{aknowledgments}{\par\addvspace{4mm}\noindent
{\it Aknowledgments. }\rm}{\par\vspace{1mm}}

\DeclareMathOperator{\tr}{Tr}

\begin{document}

\begin{frontmatter}



\title{Construction of self-dual normal bases\\
       and their complexity}

\thanks[1]{Part of this work was completed when the author was visiting the University of Limoges, funded by the London Mathematical Society Cecil King Travel Scholarship}
\author[label1]{Fran\c{c}ois Arnault}\ead{francois.arnault@unilim.fr}, \author[label2]{Erik Jarl Pickett\corauthref{1}}\ead{erikjarl.pickett@epfl.ch}, \author[label1]{St\'ephane Vinatier}\ead{stephane.vinatier@unilim.fr}

\bibliographystyle{elsart-num-sort}
\address[label1]{XLIM UMR 6172 CNRS - Universit\'e de Limoges, 123 avenue Albert Thomas, 87060 Limoges cedex, France}
\address[label2]{Math\'ematiques, \'Ecole Polytechnique F\'ed\'erale de Lausanne, 1015 Lausanne, Switzerland}

\begin{abstract}
Recent work of Pickett has given a construction of self-dual normal
bases for extensions of finite fields, whenever they exist. In this
article we present these results in an explicit and constructive
manner and apply them, through computer search, to identify the lowest
complexity of self-dual normal bases for extensions of low degree. 
Comparisons to similar searches amongst normal bases show that the
lowest complexity is often achieved from a self-dual normal basis.  
\end{abstract}

\end{frontmatter}

\section*{Introduction}
Let $q$ be a power of a prime, $n$ an integer, and let $\F_{q}$ be
the field of $q$ elements. The Galois group $G$ of the
extension $\F_{q^n}/\F_q$ is a cyclic group, generated by the
Frobenius automorphism $\phi:x\mapsto x^q$. 

A basis for $\F_{q^n}/\F_q$ consisting of the orbit of a single
element $\alpha$ under the action of the Frobenius is known as a normal
basis. In such a basis, exponentiation by $q$ is a cyclic shift of
coordinates, hence is straightforward as well as trace 
computation. The difficulty of multiplying two elements written as linear
combinations of the conjugates of $\alpha$ is measured by the
so-called \textit{complexity} of $\alpha$, defined as the number of
non zero entries in the multiplication-by-$\alpha$ matrix
\cite[\S4.1]{Menezes}. It has been 
shown in \cite{Mullin_et_al} to be at least $2n-1$, in which case the
basis is called \textit{optimal}, but this occurs only for very special
values of $n$ \cite{Gao+Lenstra}. 

The search for normal bases with low complexity has taken two
complementary directions. On the theoretical side, several authors
have attempted to build them either from roots of unity in 
larger extensions, using Gauss periods
\cite{Ash_Blake_Vanstone,Gao+Lenstra,Liao_Feng} or traces of
optimal normal bases \cite{Christopoulou_et_al}, again with some
limitations on the degree; or from the extension itself, using 
division points of a torus \cite{Blake_Gao_Mullin,Gao} or of an elliptic curve
\cite{Couveignes_Lercier}. In the latter case the authors show that fast
arithmetic can be implemented using their bases, as was also shown to
be the case for normal bases generated by Gauss periods in
\cite{Gao_et_al}.

On the experimental side, exhaustive searches of all normal bases of a
given extension have been carried out. Mullin, Onyszchuk, Vanstone and
Wilson have given a first list of lowest complexities in degree less
than $30$ over $\F_2$ in \cite{Mullin_et_al}; this list was
extended up to degree $33$ by Geiselmann \cite[Table
  5.1]{Jungnickel1993}. In odd characteristic, Blake, Gao and Mullin
computed the lowest complexities of normal bases for a handful of
small degree extensions \cite{Blake_Gao_Mullin}.  
Recently, Masuda, Moura, Panario and Thomson have reached
degree $39$ over $\F_2$ and given appealing statistics and conjectures
about the distribution of complexities \cite{Masuda_et_al}. These
authors point out that the cost of the exhaustive enumeration of the
elements of $\F_{2^n}$ used to look for normal basis generators is a
severe limitation to their method when the degree grows. On the other
hand, their Table 4 shows that the minimal complexity for normal bases
is very often reached by so-called self-dual bases (in all degrees not
divisible by $4$ up to $35$ apart from $7$, $10$, $21$). Restricting
to self-dual normal bases enables one to push computations further;
Geiselmann was indeed able to compute the lowest complexity for self-dual
normal bases over $\F_2$ up to degree $47$
\cite[\textit{loc. cit.}]{Jungnickel1993}. Comparing his results and
\cite[Table 5]{Masuda_et_al}, we see that the best found complexity for
normal bases in degree over $40$, obtained by theoretical constructions or random search,
is also reached by a self-dual normal basis for odd degrees up to $47$.


A normal basis $(\alpha,\alpha^q,\ldots,\alpha^{q^{n-1}})$ for the
extension~$\F_{q^n}/\F_q$ is said to be self-dual if
$\tr(\alpha^{q^i}\alpha^{q^j})=\delta_{i,j}$ for $0\leq i,j\leq n-1$,
where $\tr$ is the trace map from $\F_{q^n}$ to $\F_q$ and $\delta$ is
the Kronecker delta; a self-dual basis is indeed equal to its dual basis
(see \cite[\S1.2]{Menezes} for a definition), and its complexity is
the number of non zero entries in the matrix:
$$\big(\tr(\alpha\alpha^{q^i}\alpha^{q^j})\big)_{0\le i,j\le n-1}\enspace.$$
Self-dual normal bases are useful for arithmetic and Fourier
transform, and have applications in coding theory and
cryptography. Contrary to normal bases, not all extensions of finite
fields admit self-dual normal bases, but the existence conditions,
recalled in Theorem \ref{SDNBexistence} below, are mild. 
The theoretical techniques used to
construct normal bases with low complexity sometimes yield self-dual
normal bases, see \textit{e.g.} \cite[\S5.4]{Gao} or
\cite[\S5]{Blake_Gao_Mullin}, \cite[Corollary 3.5]{Gao_et_al},
\cite[Theorem 5]{Christopoulou_et_al}, \cite{Nogami_et_al}. 


In this paper we focus on the experimental side
and give the lowest complexity of self-dual normal bases in various
characteristics and degrees. 
At present, the only known strategy to
reach this goal is to compute the complexity of all the self-dual
normal bases of the extension (unless it admits an optimal self-dual
normal basis, which is easily predictable using \cite[Theorem
  2]{Liao_Sun}). In order to do so, we first construct a self-dual
normal basis for the extension, then act on it by the
\textit{orthogonal circulant} group, namely the group of change of
self-dual normal basis matrices. This group has been extensively
studied, with accurate descriptions being given in
\cite{ByrdVaughan,J-B-G,MacWilliams}. Its size is in 
$O(q^{{n}/{2}})$ (see Remark \ref{size} below), roughly the square
root of the number of normal bases in view of \cite[Corollary
  4.14]{Menezes}. It follows that exhaustive enumeration of self-dual 
normal bases is easier than that of normal bases. We shall restrict
ourselves to extensions $\F_{q^n}/\F_q$ which are either
\textit{semi-simple} (the degree $n$ prime to the characteristic $p$)
or \textit{ramified} ($n$ a power of $p$), the description of the
orthogonal circulant group in the ``mixed'' case being a bit more elaborate. 

We now describe our work more precisely. First we recall the 
necessary and sufficient conditions for the existence of self-dual normal bases
\cite{LempelWeinberger}:  
\begin{ithm}[Lempel-Weinberger]
\label{SDNBexistence}
The extension of finite fields $\F_{q^n}/\F_q$ has a self-dual normal
basis if and only if 
either the degree~$n$ is odd or $n\equiv2$ modulo~4 and $q$ is even.
\end{ithm}
The existence proof in \cite{LempelWeinberger} is constructive in the
sense that, given a normal basis for the extension, it
describes a procedure to transform it into a self-dual normal basis.
Wang proposed another transformation procedure in \cite{Wang} when
$q=2$ and $n$ is odd, involving solving a system of equations. Poli
extended Wang's method to deal with the general characteristic $2$
case in \cite{Poli}. Recently, Pickett designed in \cite{Pickett2} a
construction that extends the former ones to the odd characteristic
case, dealing separately with the semi-simple case and the ramified case.

The construction of a normal basis for a given extension is well known
and widely implemented. Therefore, the methods described above enable one to
construct a self-dual normal basis under the existence conditions of
Theorem \ref{SDNBexistence}. To
our knowledge, this has not been implemented before, except in the restrictive
case in which Wang's method applies. In this paper we apply
Pickett's construction to compute a self-dual normal basis of a
given extension whenever it exists. Note that for this first goal the
method in \cite{LempelWeinberger} is simpler and faster, but most of
the computations involved in Pickett's construction must be
implemented if one wants to compute the action of the
orthogonal circulant group as well.

The criterion used in \cite{Wang} to determine which changes of basis are
appropriate has been generalised to any characteristic and degree, see
\cite[Lemma 5.5.3]{Jungnickel1993}, where it is expressed in terms of
circulant  matrices. Here we restate it in terms of the group 
algebra~$\F_q[G]$ as in \cite{Pickett2}.
Conjugation $u\mapsto\overline{u}$ in~$\F_q[G]$ is the~$\F_q$-algebra
automorphism obtained from $g\mapsto g^{-1}$ for all $g\in G$; if
$u=\sum_{k=0}^{n-1}u_k\phi^k\in\F_q[G]$ and $\alpha\in\F_{q^n}$, we put
$u\circ\alpha=\sum_{k=0}^{n-1}u_k\phi^k(\alpha)\in\F_{q^n}$.
\begin{ithm}\label{vv}
Assume that $\alpha$ is a generator of a normal basis 
of~$\F_{q^n}$ over~$\F_q$ and let
$$R = \sum_{g\in G} \tr\big( \alpha g(\alpha)\big)\,g\
  \in\ \F_q[G]\enspace.$$
Any $v\in\F_q[G]$ such that $v\overline{v}=R$ is invertible, and the
map $v\mapsto v^{-1}\circ\alpha$ is a one-to-one correspondence
between the set of solutions of the equation $v\overline{v}=R$ in
$\F_q$ and the set of elements of $\F_{q^n}$ that generate a self-dual
normal basis. 
\end{ithm}
In Section \ref{groupalgebra} we first explain how this result can be
deduced from the statement on circulant matrices
\cite[\textit{loc. cit.}]{Jungnickel1993}. Our main interest is in
implementing Pickett's method as an algorithm, and since the language
he uses to describe his construction of a solution of the equation
$v\overline v=R$ in \cite[\S3]{Pickett2} is quite elaborate --- his
framework is wider than ours --- we reformulate it in terms of the
polynomial ring $\F_q[X]/(X^n-1)$; the resulting algorithm to compute
a self-dual normal basis is described in the last section. We
remark that this construction gives an alternative proof of the
sufficiency of the conditions of Theorem~\ref{SDNBexistence}; for
interest we give a proof of their necessity, mainly based on
Theorem \ref{vv}, and simpler than the original (see
\cite[Propositions 4.3.4 and 5.2.2]{Jungnickel1993}). 

Section \ref{change} deals with the orthogonal
circulant group $O(n,q)$. Its elements are the $n\times n$ matrices
$P$ over $\F_q$ that are circulant ($P_{i+k\!\!\mod n,j+k\!\!\mod n}=P_{i,j}$
for $0\le i,j,k\le n-1$) and orthogonal ($P^t\cdot P=I$, where $P^t$
is the transpose matrix of $P$ and $I$ the identity $n\times n$
matrix). It follows from Theorem \ref{vv} that $O(n,q)$ is isomorphic to
the subgroup of $\F_q[G]^\times$ consisting of the solutions of the
equation $v\overline v=1$. In both the semi-simple and the ramified
case we indicate how this equation can be solved; the resulting
algorithms are described in the last section. Doing so we recover the
number of self-dual normal bases, as derived in \cite{J-B-G,JMV} from
MacWilliams' results about the orthogonal circulant group
\cite{MacWilliams} (see 
\cite[5.3]{Jungnickel1993} for a summary). In the ramified (and odd
characteristic) case our construction is a variation, adjusted to our
situation, of MacWilliams' iterative construction; we also present a
new explicit formula for the solutions. 

In Section \ref{experiments} we present our algorithms, experimental
results and conclusions. For semi-simple extensions in odd
characteristic, the lowest complexity we find is close to that
obtained for normal bases from exhaustive computer search 
\cite{Blake_Gao_Mullin} or from theoretical constructions
\cite{Liao_Feng}, as this was already the case in even
characteristic. We also observe an interesting behaviour under base
field extension. 
When the extension is of degree $p$ in odd characteristic $p$ we
recover the basis with very low complexity $3p-2$ described in
\cite{Blake_Gao_Mullin}. 

\section{Construction of a self-dual normal basis}\label{groupalgebra}
Our algorithm to find a self-dual normal basis relies on the
interpretation in terms of polynomial rings of Pickett's 
construction of a solution $v$ of the equation $v\overline v=R$ of Theorem
\ref{vv} (under the necessary conditions of Theorem \ref{SDNBexistence}).
The majority of this section is devoted to presenting this
interpretation. First, however, we deduce Theorem \ref{vv} from statements in terms of circulant matrices.
At the end of the section we show how 
to deduce the necessity of the conditions of Theorem
\ref{SDNBexistence} from Theorem \ref{vv}. 
\subsection{Proof of Theorem \ref{vv}}
Consider the one-to-one correspondence between $\F_q[G]$ and circulant
$n\times n$ matrices over $\F_q$, given by
\begin{equation}\label{circulant}
v=\sum_{j=0}^{n-1}\rho_j\phi^j\in\F_q[G]\mapsto
C_v=\big(\rho_{j-i\bmod n}\big)_{0\leq i,j\leq n-1}\enspace.
\end{equation}
One has $C_1=I$ and, for any $v,w\in\F_q[G]$, $C_v\cdot C_w=C_{vw}$,
so (\ref{circulant}) yields a group isomorphism between
$\F_q[G]^\times$ and the abelian group of invertible circulant
$n\times n$ matrices over $\F_q$. Note that the matrix 
$C_R=\big(\tr(\alpha^{q^i+q^j})\big)$ is invertible since $\alpha$
generates a normal basis, see \cite[Corollary 1.3]{Menezes}, so
$R\in\F_q[G]^\times$ and $v\overline v=R$ implies $v$ invertible as
well. 

Moreover one has $C_{\overline v}=(C_v)^t$, where $(C_v)^t$ is the
transpose matrix of $C_v$. It follows that the equation
$v\overline v=R$ is equivalent to  
\begin{equation}\label{selfdualmatrix}
  C_v\cdot(C_v)^t=\big(\tr(\alpha^{q^i+q^j})\big)_{0\leq i,j\leq n-1}\enspace.
\end{equation}

For $x\in\F_{q^n}$, let $[x]$ denote the $n\times n$ matrix whose
$j$-th column, $0\le j\le n-1$, consists of the coordinates of
$x^{q^{j}}$ in a fixed $\F_q$-basis of $\F_{q^n}$. Then one has, for
any $v\in\F_q[G]$, $x\in\F_{q^n}$: 
$$[v\circ x]=[x]\cdot C_v\enspace.$$

Let $P$ be some invertible $n\times n$ matrix over $\F_q$, then
the columns of $B=[\alpha]P$ are the coordinates in the fixed
$\F_q$-basis of $\F_{q^n}$ of a normal basis if and only if $P$ is a
circulant matrix, see \cite[Lemma 3.1.3]{Jungnickel1993}. Further, for
such a $P$, its inverse $P^{-1}$ is also circulant and from
\cite[Lemma 5.5.3]{Jungnickel1993} we know that the columns of $B$ form a 
self-dual normal basis if and only if 
\begin{equation}\label{selfdualmatrix2}
P^{-1}\cdot(P^{-1})^t
  =
  \big(\tr(\alpha^{q^i+q^j})\big)_{0\leq i,j\leq n-1}\enspace.
\end{equation}

 If $v\overline v=R$, then
$C_v$ is circulant invertible and 
$(C_v)^{-1}=C_{v^{-1}}$ satisfies (\ref{selfdualmatrix2}), hence
$B=[\alpha]C_{v^{-1}}=[v^{-1}\circ\alpha]$ is a self-dual normal
basis; if $\beta$ generates a self-dual normal basis, let $P$ be such
that $[\beta]=[\alpha]P$, it is circulant and so is its inverse, and
by (\ref{selfdualmatrix}) the element
$v\in\F_q[G]$ such that $P^{-1}=C_v$ satisfies $v\overline v=R$. These
two maps are clearly mutual inverses, which completes the proof.

\subsection{Interpretation of Pickett's construction in terms of polynomial rings}

The Galois group $G$ of $\F_{q^n}$ over~$\F_q$ is
cyclic of order~$n$ and generated by the Frobenius $\phi$, so we may
identify the $\F_q$-algebras $\F_q[G]$ and  $\F_q[X]/(X^n-1)$ through
the isomorphism mapping $\phi$ to $X$.

Write $n=p^en_1$, where $p$ is the characteristic of~$\F_q$ and $n_1$
is prime to $p$. We take advantage of the following result 
\cite[Theorems 3.3.13 and 5.1.9]{Jungnickel1993} to split the
extension into two parts.  
\begin{lemma}\label{compositum}
Let $m,n$ be two co-prime integers. Suppose $\alpha$ (resp. $\beta$) is
a generator of a self-dual normal basis of $\F_{q^m}$
(resp. $\F_{q^n}$) over $\F_q$, then $\alpha\beta$ is a generator of a
self-dual normal basis of the compositum $\F_{q^{mn}}$ over
$\F_q$. Moreover, the complexity of $\alpha\beta$ is the product of
the complexities of $\alpha$ and of $\beta$.
\end{lemma}
By the former result,
we may deal separately
with the two cases $n=p^e$ which we call the ramified case, and $n$
co-prime to $p$, the so-called semi-simple case. 
We show how to construct a solution $v$ of the equation $v\overline v=R$ of
Theorem \ref{vv} in each of
these two cases, under the existence conditions of a self-dual normal
basis of Theorem \ref{SDNBexistence}. 
Multiplying the bases obtained this way then yields self-dual
normal bases for the extensions with ``mixed degree'' $n=n_1p^e$ with
$n_1\ge 2$ and $e\ge 1$. 
\subsubsection{The ramified case ($n=p^e$)}\label{ram_section}

In this case, the algebra~$\F_q[G]$ is isomorphic to
$\F_q[X]/(X-1)^n$.  Let $\epsilon:\F_q[G]\rightarrow\F_q$ be the 
augmentation map given by
$\epsilon(\sum_{k=0}^{n-1}a_k\phi^k)=\sum_{k=0}^{n-1}a_k$. This is a homomorphism of
$\F_q$-algebras whose kernel is a codimension $1$ subspace of
$\F_q[G]$. Further $\epsilon(\sum_{k=0}^{n-1}a_k\phi^k)=0$ implies
$\sum_{k=0}^{n-1}a_k\phi^k=\sum_{k=0}^{n-1}a_k(\phi^k-1)$, and
therefore the kernel is $(\phi-1)\F_q[G]$. Invertible 
elements in $\F_q[G]$ are those which have non-zero image 
under the map~$\epsilon$ (because invertible modulo~$(X-1)^n$ means
invertible modulo~$X-1$), hence the group 
$\F_q[G]^\times$ has order~$q^{n-1}(q-1)$.  In fact, it is the direct
product of $\F_q^\times$ by $U=1+(\phi-1)\F_q[G]$, the inverse image of~$1$ under the map~$\epsilon$.

Under the necessary conditions of Theorem \ref{SDNBexistence}, we have
two cases to consider.
\begin{proposition}\label{SDNBramified}
Recall that $p$ is the characteristic of $\F_q$.
If $p=n=2$, $\beta\in\F_{q^2}$ generates a self-dual normal basis if
and only if $\tr(\beta)=1$.
If $p$ is odd and $n=p^e$, there exists $\omega\in\F_q[G]$ such that
$\omega^2=R$; further one then has $\omega=\overline\omega$.
\end{proposition}
\begin{proof}
The even characteristic case is straightforward. 
We proceed with the odd characteristic case.
Recall that $R\in\F_q[G]^\times$ and note that $\overline{R}=R$. One can easily see that
$\epsilon(R)=\tr(\alpha)^2$ (detailed in the proof of Lemma \ref{Erik3.5}
below), so that the decomposition of~$R$ in the 
above direct product is $R=\tr(\alpha)^2\cdot(1+(\phi-1)R')$ for some 
$R'\in\F_q[G]$.  The second factor is also a square as it belongs to
the group~$U$ which is of odd order, hence
$R=\omega^2$ for some~$\omega$. Further $\overline R=R$ implies 
$\overline{\omega}^2=\omega^2$, so that $\overline\omega/\omega$ is a
square root of~$1$ living in the group~$U$ of odd order.  Thus
$\overline{\omega}=\omega$. 
\end{proof}

\subsubsection{The semi-simple case ($\gcd(n,q)=1$)} \label{ss_section}

We assume that $n$ is odd to fit with the conditions of
Theorem~\ref{SDNBexistence} (but $q$ could be odd or even). The
polynomial $X^n-1$ is square free and has monic irreducible factors
over~$\F_q$~: 
\begin{equation}\label{polydecomposition}
X^n-1 = \prod_{i=1}^\sigma f_i(X) \prod_{j=1}^\tau g_j(X) \cdot g_j^*(X)
\end{equation}
where $g_j^*$ denotes the reciprocal polynomial (up to a constant) of~$g_j$ and where the
$f_i$ are the self-reciprocal (also up to a constant) irreducible factors.  We will now express the
equation $R=v\overline v$ in this decomposition, solve it, and then lift back the solution
to~$\F_q[G]$.

Let $m$ be the order of $q$ modulo~$n$.  The field~$\F_{q^m}$ contains
a primitive $n$-th root~$\zeta$ of~1.  On the set $\{0,\ldots,n-1\}$
we define the {\it cyclotomic equivalence relation:} $s\sim s'$ if
there exists $k$ such that $s\equiv q^ks'$ mod~$n$.
Note that $0$ forms a class on its own and that the integers prime
to $n$ belong to classes with the same cardinal equal to the order of
$q$ modulo $n$. Namely, since $n$ and $q$ are co-prime, the cyclotomic
equivalence relation restricts to $(\Z/n\Z)^\times$ and for $s,s'$
invertible modulo $n$, $s\sim s'$ if and only if $s$ and $s'$ belong
to the same coset in $(\Z/n\Z)^\times/\langle q\rangle$.

The following proposition justifies the terminology.

\begin{proposition}\label{cyclotomicfactorization} 
(a) If $\zeta^s$ is a root of an irreducible factor of~$X^n-1$, then the
other roots are the $\zeta^{s'}$ where $s'\sim s$.\\
\noindent(b) The $\zeta^s$ such that $s\sim(n-s)$ are roots of a
self-reciprocal factor $f_i$. 
The $\zeta^s$ such that $s\not\sim n-s$ are roots of a non
self-reciprocal factor $g_j$.\\ 
\noindent(c) The number of cyclotomic classes is equal to the number $\sigma+2\tau$ of irreducible
factors of~$X^n-1$.\\
\noindent(d) The self-reciprocal factors $f_i$ have even degree,
except $f_1=X-1$.
\end{proposition}

\begin{proof} 
(a), (b), (c) are clear.  Let us prove~(d).  If $\zeta^s$ is a root of
  an~$f_i$, then $\zeta^{n-s}$ is also a root.  If we exclude the case $s=0$
  corresponding to the factor~$X-1$, the two roots 
$\zeta^s$ and~$\zeta^{n-s}$ are distinct, because $n$ is odd.  Hence
  $f_i$ has en even number of roots in an algebraic closure.  
\end{proof}

From the Chinese Remainder Theorem, the algebra $\F_q[X]/(X^n-1)$ is
isomorphic to a product of $\sigma+2\tau$ fields:
\begin{equation}\label{CRT}
  \F_q[X]/(X^n - 1)
  \simeq
  \prod_{i=1}^\sigma \F_q[X]\big/\big(f_i(X)\big)
  \times
  \prod_{j=1}^\tau \Big(\F_q[X]\big/\big(g_j(X)\big) \times \F_q[X]\big/\big(g_j^*(X)\big)\Big)\enspace.
\end{equation}
Each factor in the RHS of this equation is a an extension of~$\F_q$
contained in~$\F_{q^m}$ (recall $m$ is the order of $q$ modulo
$n$). The evaluation map $u(X)\in\F_q[X]/(f)\mapsto
u(\zeta^s)\in\F_q(\zeta^s)$, where $f$ is an $f_i$ or a $g_j$ and
$s\in\{0,\ldots n-1\}$ such that $f(\zeta^s)=0$, is a field
isomorphism. We obtain the following result:

\begin{proposition}
Let $S$ be a set of representatives of cyclotomic classes.  The map
\begin{equation}\label{Adecomposition}  
  \left\{\begin{array}{rl}
     \F_q[X]/(X^n-1) &\longrightarrow \prod_{s\in S} \F_q(\zeta^s) \\
     u(X)            &\longmapsto \big(u(\zeta^s)\big)_{s\in S} 
\end{array}\right.
\end{equation}
is an $\F_q$-algebra isomorphism.
\end{proposition}  
For practical reasons (mainly to deal with square matrices), we also 
consider the map~$\FT$ (a Fourier Transform)
\begin{equation}\label{FT}
  \FT : 
  \left\{\begin{array}{rl}
     \F_q[X]/(X^n-1) &\longrightarrow (\F_{q^m})^n  \\
     u(X)            &\longmapsto \big(u(\zeta^s)\big)_{0\leq s\leq n-1}  \end{array}
  \right.
\end{equation}
which is a homomorphism of $\F_q$-algebras, with matrix
$F(\zeta)=(\zeta^{ij})_{0\leq i,j\leq n-1}$. Compared with isomorphism
(\ref{Adecomposition}), we now compute a component at every $0\le s\le
n-1$; the components corresponding to indices in the same coset under
$\sim$ are cyclically permuted when applying the Frobenius $\phi$.

We note the following easy but useful relation involving the matrices
$F(\zeta)$ and $F(\zeta^{-1})=(\zeta^{-ji})_{0\leq i,j\leq n-1}$: 

\begin{lemma}  
$F(\zeta^{-1})F(\zeta)=nI$.
\end{lemma} 
As a consequence, the following linear map~$\FTI$, with matrix $F(\zeta^{-1})$, can be used to
compute the inverse of $\FT$.
\begin{equation}\label{FTI}
  \FTI : 
  \left\{\begin{array}{rl}
     (\F_{q^m})^n       &\longrightarrow \F_{q^m}[X]/(X^n-1)  \\
     (r_0,\ldots,r_{n-1})  &\longmapsto \sum_{t=0}^{n-1} u_tX^t 
                \text{\ \ \ \ where $u_t=\sum_{i=0}^{n-1}
                  r_i\zeta^{-ti}\ $.} \end{array}
  \right.
\end{equation}
This is because $\FTI\big(\FT(u)\big)=nu$ for each $u\in\F_q[X]/(X^n-1)$.


The idea here is to express $R$ as an element of the RHS
of~(\ref{Adecomposition}), to solve the equation in each component,
and to bring back the solution to~$\F_q[X]/(X^n-1)$.  The
conjugation map, induced by $X\mapsto X^{n-1}$ in~$\F_q[X]/(X^n-1)$
is given by $\zeta\mapsto\zeta^{-1}$ and will sometimes be 
denoted by~$J$ in the RHS of~(\ref{Adecomposition}).

Let $R$ be as in~Theorem \ref{vv}. The $s$-coordinate of~$\FT(R)$ is
$R_s=\sum\limits_{i=0}^{n-1}\tr(\alpha^{1+q^i})\zeta^{si}$. 

We begin with the cyclotomic class $s=0$.
Here, $\F_q(\zeta^s)=\F_q$ and the conjugation map~$J$ acts
trivially. Note that $R_0=\epsilon(R)$.

\begin{lemma} [3.5 in~\cite{Pickett2}] \label{Erik3.5}
With $v_{0}=\tr(\alpha)$, we have 
$v_{0}\overline{v_{0}}=R_{0}$. 
\end{lemma} 

\begin{proof} 
We have $J\big(\tr(\alpha)\big)=\tr(\alpha)$ and 
$$
  \tr(\alpha)^2 
  =
  \Big(\sum_{i=0}^{n-1} \alpha^{q^i}\Big)^2
  =
  \sum_{i,j=0}^{n-1} \alpha^{q^i+q^j}
  =
  \sum_{i,k=0}^{n-1} \alpha^{q^i(1+q^k)}
  =
  \sum_{k=0}^{n-1} \tr(\alpha^{1+q^k})
  =
  R_{0}.
$$\end{proof}

We now consider the cyclotomic classes $s$ such that $s\not\sim n-s$.
\begin{lemma}[3.6 in~\cite{Pickett2}]\label{Erik3.6}  
Let $s'\in S$ such that $s'\sim n-s$.   We have $R_s=R_{s'}$.
Putting $v_{s,s'}=(R_s,1)\in\F_q(\zeta^s)\times\F_q(\zeta^{s'})$, we have 
$v_{s,s'}\overline{v_{s,s'}}=(R_s,R_{s'})$.
\end{lemma}

\begin{proof} 
The conjugation map~$J$ exchanges coordinates in $\F_q(\zeta^s)\times\F_q(\zeta^{s'})$: 
$J(u,u^*)=(u^*,u)$.  As $R$ is invariant by conjugation, we have $R_s=R_{s'}$.  Therefore 
$v_{s,s'}J(v_{s,s'})=(R_s,1)(1,R_s)=(R_s,R_{s'})$.  
\end{proof}

We finally deal with the cyclotomic classes $s$ such that $s\not=0$ and $s\sim n-s$.
\begin{lemma}[3.7 in~\cite{Pickett2}]\label{Erik3.7}
Let $s\in S$ such that $0\neq s$ and $s\sim n-s$.
Then the field $\F_q(\zeta^s)$ is stable under the conjugation
map~$J$, and we denote by $\F_q(\zeta^s)^J$ the fixed subfield. Further
$R_s$ (resp. $-R_s$) has a square root $u$ 
(resp. $u'$) in $\F_q(\zeta^s)$.
We consider three cases:
\begin{enumerate}[(a)]
\item the case where $u\in\F_q(\zeta^s)^J$, then $v_s=u$ satisfies
$v_s\overline{v_s}=R_s$;
\item the case where $u'\notin\F_q(\zeta^s)^J$, then $v_s=u'$ satisfies
$v_s\overline{v_s}=R_s$;
\item the case where $u\notin\F_q(\zeta^s)^J$ and $u'\in\F_q(\zeta^s)^J$, then there
exists an integer $n$ such that $-n$ is a non-zero square $\eta^2$
modulo the characteristic $p$ of~$\F_q$, but $-(n-1)$ is not a square
modulo~$p$, and there exists an integer $\nu$ such that $\nu^2\equiv n-1$ 
modulo~$p$.  We put $v_s=(\nu u+u')/\eta$, then
$v_s\overline{v_s}=R_s$.
\end{enumerate}
\end{lemma}

\begin{proof} 
From Proposition~\ref{cyclotomicfactorization}, the field $\F_q(\zeta^s)$ is some
extension~$\F_{q^r}$ over~$\F_q$ with $r$ even.  We have $\overline{\zeta^s}=\zeta^{-s}\neq\zeta^s$
because $n$ is odd.  Hence $J$ restricted to $\F_q(\zeta^s)$ is an order
$2$ field automorphism, which by Galois theory defines a unique index
$2$ subextension $\F_q(\zeta^s)^J=\F_{q^{r/2}}$.
Note that, each element of~$\F_q(\zeta^s)^J$ is a square
in~$\F_q(\zeta^s)$ (because $(q^r-1)/(q^{r/2}-1)=q^{r/2}+1$ is even).  
Both $R_s$ and $-R_s$ are invariant under $J$, hence
they are both squares in~$\F_q(\zeta^s)$.  

If $\overline u=u$, namely in case (a), then $u\overline u=u^2=R_s$;
if $\overline{u'}\not=u'$, namely in case (b), then
$\overline{u'}=-u'$ and $u'\overline{u'}=-u'^2=R_s$. 

Suppose now (case
c) that $\overline u=-u$ and $\overline{u'}=u'$.  As $-1=-R_s/R_s$, we
know that $-1$ is not a square in~$\F_q(\zeta^s)^J$, nor in~$\F_p$.
Hence the first $n>1$ such that $-n$ is a square modulo~$p$ 
exists and satisfies the required conditions.  Also, because neither
$-1$ nor $-(n-1)$ are squares modulo~$p$, there exists an integer $\nu$
such that $\nu^2\equiv(n-1)$ modulo~$p$. Taking the residues of $\eta$
and $\nu$ modulo $p$, we have $\overline\eta=\eta$ and
$\overline\nu=\nu$ because $\F_p\subseteq \F_q(\zeta^s)^J$.   With
$v_s=(\nu u+u')/\eta$, we have 
$\overline{v_s}=(-\nu u+u')/\eta$ and it follows that
$v_s\overline{v_s}=(-\nu^2u^2+u'^2)/\eta^2=(-(n-1)R_s-R_s)/(-n)=R_s$.  
\end{proof}

We have solved the equation $v_s\overline{v_s}=R_s$ for every
cyclotomic class $s$, thus by the $\F_q$-algebra isomorphism
(\ref{Adecomposition}) we get a solution $v\in\F_q[G]$ of the equation 
$v\overline v=R$. 
\subsection{The necessity of the conditions of Theorem \ref{SDNBexistence}}

If $\alpha$ is a generator of a self-dual normal basis of
$\F_{q^{nm}}$ over~$\F_q$, then $\tr_{\F_{q^{nm}}/\F_{q^{n}}}(\alpha)$
is a generator of a self-dual normal basis of~$\F_{q^{n}}$ over~$\F_{q}$, see \cite[Lemma 4.3]{Pickett2}. Therefore, to prove the necessity of the conditions in Theorem \ref{SDNBexistence} we need just consider the cases $\F_{q^2}/\F_q$ for $q$ odd and $\F_{q^4}/\F_q$ for $q$ even.


When $q$ is odd, $\tr(\alpha\alpha^q)=2N(\alpha)$ for any $\alpha\in\F_{q^2}$, where
$N(\alpha)$ denotes the norm of $\alpha$ in the extension, hence 
$\tr(\alpha\alpha^q)=0$ would imply $\alpha=0$.


Let $q$ be even, and assume for contradiction that there exists a normal
basis generator $\alpha$ of $\F_{q^4}/\F_q$ and an element $v\in\F_q[G]$ such that $v\overline 
v=\tr(\alpha^2)+\tr(\alpha\alpha^q)\phi+\tr(\alpha\alpha^{q^2})\phi^2+\tr(\alpha\alpha^{q^3})\phi^3$. Note
that $\tr(\alpha\alpha^{q^3})=\tr(\alpha\alpha^q)$ and
$\tr(\alpha\alpha^{q^2})=2\tr_{\F_{q^2}/\F_q}\big(N_{\F_{q^4}/\F_{q^2}}(\alpha)\big)=0$. Writing
$v=a+b\phi+c\phi^2+d\phi^3$ with $a,b,c,d\in\F_q$ and letting
$\beta=\alpha+\alpha^{q^2}$, we easily get the equations:
$$a+b+c+d=\tr(\alpha)=\beta+\beta^q\quad,\quad(a+c)(b+d)=\tr(\alpha\alpha^q)=\beta\beta^q\enspace.$$
It follows that $\{\beta,\beta^q\}=\{a+c,b+d\}$, namely
$\beta\in\F_q$, which is impossible since it
would imply $\alpha+\alpha^{q^2}=\alpha^q+\alpha^{q^3}$, contradicting
the fact that $\alpha$ generates a normal basis. The result now
follows using Theorem \ref{vv}.

\section{Change of self-dual normal basis}\label{change}

The next result, which is essentially a different formulation of the
``key'' lemmas 2 and 3 of \cite{JMV}, is an immediate consequence of
Theorem \ref{vv} and the observations that if $\alpha$ generates a
self-dual normal basis, then $R=1$, and that if $v\overline{v}=1$, then
$v^{-1}=\overline{v}$. 
\begin{corollary}\label{sdnbvp1}
Let $\alpha$ generate a self-dual normal basis of~$\F_{q^n}$
over~$\F_q$. The map $v\mapsto\overline v\circ\alpha$ is an
isomorphism between the group of solutions of the equation $v\overline
v=1$ in $\F_q[G]$ and the group of elements of $\F_{q^n}$ that
generate a self-dual normal basis.
\end{corollary}

It follows that computing all self-dual normal bases from
one is equivalent to finding all the solutions $v\in\F_q[G]^\times$ of
the equation $v\overline v=1$. We devote the rest of this section to
explain how this equation can be solved, first in the semi-simple case
and then in the ramified case. 
\subsection{The semi-simple case}

The decomposition~(\ref{Adecomposition}) from Section
\ref{groupalgebra} is useful to find the solutions of this equation.
Let $V(X)\in\F_q[X]/(X^n-1)$.   

\begin{proposition}\label{ssc}
The polynomial $V(X)$ satisfies the equation $V(X)V(X^{n-1})=1$ modulo~$X^n-1$
if and only if the following conditions hold:  
$$
  \left\{\begin{array}{ll}
    V(1) = \pm1          &\text{\ \ (case $s=0$),} \cr
    V(\zeta^s)V(\zeta^{-s})=1 &\text{\ \ for $s\not\sim n-s$,}\cr
    V(\zeta^{s})^{q^{r/2}+1}=1     &\text{\ \ for $0\neq s\sim n-s$, 
              where $r$ is such that $\F_q(\zeta^s)=\F_{q^r}$}.\cr
\end{array}\right.
$$
Note that $r$ is the degree of the irreducible factor $f_i$ of~$X^n-1$
such that $f_i(\zeta^s)=0$. 
\end{proposition}

\begin{proof} 
The component at $s=0$ is $V(1)$ and the equation we need to solve in
$\F_q(\zeta^0)=\F_q$ is simply $V(1)^2=1$ because the action of
conjugation in $\F_q$ is trivial. 

For $s\not\sim n-s$, we have to consider the product $\F_q(\zeta^s)\times\F_q(\zeta^{-s})$.  We
have seen in the proof of Lemma~\ref{Erik3.6} that conjugation swaps coordinates in
these two factors.  The solutions 
are the powers of $(g_s,g_s^{-1})$ were $g_s$ is any primitive element of the~$\F_q(\zeta^s)$.

For $0\neq s\sim n-s$, we have seen in the proof of
Lemma~\ref{Erik3.7} that the set of invariants under 
conjugation~$J$ is the subfield~$\F_{q^{r/2}}$ of~$\F_{q^r}=\F_q(\zeta^s)$.   Conjugation~$J$ is an
$\F_{q^{r/2}}$-automorphism of~$\F_{q^r}$ of order~2, hence $J(x)=x^{q^{r/2}}$ for $x\in\F_{q^r}$.
The equation we want to solve can be written $x^{q^{r/2}+1}=1$.  Note that $q^{r/2}+1$ divides
$q^r-1$ so we find exactly $q^{r/2}+1$ solutions, generated by any element of order~$q^{r/2}+1$
in~$\F_q(\zeta^s)$.  
\end{proof}

We remark that this proof provides generators for the group of
solutions of $v\overline v=1$, so we can easily derive the cardinality of
this group, which by Corollary~\ref{sdnbvp1} is also the number of
self-dual normal bases of $\F_{q^n}$ over~$\F_q$. As expected, this calculation
agrees with the result in~\cite{JMV} which was obtained using the
formulas given in~\cite{MacWilliams} --- note that the cyclic shift of
a basis is considered to be the same basis in~\cite{JMV}, but not
here, so our formula differs from the one found there by a factor~$n$.
\begin{theorem}\label{SDNBcount}
Consider the decomposition~(\ref{polydecomposition}) of~$X^n-1$ over~$\F_q$.  The
number of self-dual normal bases of~$\F_{q^n}$ over~$\F_q$ is given by
$$
  2^a \prod_{i=2}^\sigma (q^{c_i} + 1) \prod_{j=1}^\tau (q^{d_j} - 1)
  \text{\ \ with }
  \left\{\begin{array}{l}
       \hbox{$a=0$ for even~$q$ and $a=1$ for odd~$q$},\cr
       \hbox{$2c_i=\deg f_i$ and $d_j=\deg g_j$}.\cr
\end{array}\right. 
$$
\end{theorem}

\begin{proof} 
The case $s=0$ has solutions $\pm1$ in odd characteristic, and only~1
for even~$q$. For the case $0\neq s\sim n-s$, we found a generator of
order~$q^c+1$ for the set of solutions in the field $\F_q(\zeta)$.
For the case $s\not\sim n-s$, let $g$ be a primitive element
in~$\F_q[X]/(f)\simeq\F_q(\zeta)$, the solutions are the powers of $(g,g^{-1})$.   
\end{proof} 

\subsection{The ramified case}
We deal only with the odd characteristic case, so we let $p$ be
an odd prime number, and $q$ and $n$ be powers of $p$. 

\begin{theorem}\label{thm:MacWilliams}
There are $2q^{\frac{n-1}{2}}$ solutions $v\in\F_q[G]$ to the equation
$v\wb v=1$.
\end{theorem}
This result can easily been derived from \cite[Theorem 2]{J-B-G},
which states that if $n=sp$, where $s$ is any integer, then
$|O(sp,q)|=q^{(p-1)s/2}|O(s,q)|$. The original statement is due to
MacWilliams in the prime base field case
\cite[Theorem 2.6]{MacWilliams}. 
We now reinterpret MacWilliams' constructive proof in our specific
case: $n$ a power of $p$, so as to explain the structure of the
algorithm we used to compute the orthogonal circulant group in the
ramified case.
\begin{proof}
First note that the solutions of the equation $v\wb v=1$ all lie in
$\F_q[G]^\times$, and recall from Subsection \ref{ram_section} that
$\F_q[G]^\times$ is the direct product
$\F_q^\times\times(1+(\phi-1)\F_q[G])$, the first component being
simply the image by the augmentation map $\epsilon$. For
$v\in\F_q[G]^\times$, let $w\in(\phi-1)\F_q[G]$ be such that
$v=\epsilon(v)(1+w)$, then $v\wb v=1$ if and only if
$\epsilon(v)=\pm1$ and $w+\overline w+w\overline w=0$. Setting
$r=w+\frac{w\overline w}{2}$, the second condition becomes
$r=-\overline r$, namely 
\begin{equation}\label{r}
r=\sum_{i=1}^{\frac{n-1}{2}}r_i(\phi^i-\phi^{n-i})
\end{equation}
for some $r_i\in\F_q$, hence $r$ can take $q^{\frac{n-1}{2}}$ values
in $\F_q[G]$.
We now show that $w$ is uniquely defined by $r$, and how it can be
computed, see \cite[Appendix A]{MacWilliams}. One has
$w=-r+\frac{w\overline w}{2}$, hence $\overline w=r+\frac{w\overline
  w}{2}$ and $w\overline w=-r^2+\frac{(w\overline w)^2}{4}$, so that:
$$w=-r-\frac{r^2}{2}+\frac{(w\overline w)^2}{8}\enspace.$$
Replacing iteratively $w\overline w$ by $-r^2+\frac{(w\overline
  w)^2}{4}$ in the above formula increases the (even) power to which
$w\overline w$ appears; this process terminates since, as an element
of $(\phi-1)\F_q[G]$, $w=(\phi-1)y$ for some $y\in\F_q[G]$, so
$w^n=(\phi^n-1)y^n=0$.  
\end{proof}
\begin{remark}\label{size}
In the odd characteristic case, the formula in Theorem
\ref{SDNBcount} reads:
$$2\prod_{i=2}^\sigma (q^{c_i} + 1) \prod_{j=1}^\tau (q^{d_j} -
1)\approx 2q^{\sum_{i}{c_i}+\sum_{j}{d_j}}=2q^{(n-1)/2}\enspace.$$
In both semi-simple and ramified cases, the size of the
trace-orthogonal group is close to $2\sqrt{q^{n-1}}$, which means
that an exhaustive search quickly becomes lengthy when $q$ or $n$ increases.
\end{remark}

We now show that one can also get an explicit formula for the solutions
of the equation.
\begin{theorem}\label{thm:ramified}
The solutions $v\in\F_q[G]$ to the equation
$v\wb v=1$ are exactly the sums $v=\sum_{i=0}^{n-1}v_i(\phi-1)^i$
with $v_0=\pm 1$ and, for $1\le i\le\frac{n-1}{2}$, $v_{2i-1}$ is any
element of $\F_q$ and $v_{2i}\in\F_q$ is such that:
\begin{equation}\label{v}
\sum_{j=1}^{2i}\sum_{k=0}^{j}(-1)^k\binom{n-k}{2i-j}v_kv_{j-k}=0\enspace.
\end{equation}
\end{theorem}
Note that (\ref{v}) gives a formula for $v_{2i}$ in terms of the
$v_{k}$ with $0\le k\le 2i-1$, for instance $-2v_0v_2=-v_1^2+v_0v_1$
and $-2v_0v_4=v_0v_2-v_1v_2-2v_1v_3+v_2^2+3v_0v_3$.
Our proof begins as a specialisation to the
case $s=1$ of that of \cite[Satz 3.3]{BG} --- note that \cite{J-B-G}
points out a mistake in the end of the proof of this
statement; dealing with this simpler case enables us to 
deduce a constructive formula. 
\begin{proof}
We wish to solve the equation $v\overline v=1$ in $\F_q[G]$. We shall
proceed by successive approximation, solving $v\overline v\equiv 1$
modulo $(X-1)^i$ for $1\le i\le n$, where we identify again $v$ and its
image under the isomorphism
$$\F_q[G]\cong\F_q[X]/(X-1)^n$$
mapping $\phi$ to $X$. The
first step is obvious: $\F_q[X]/(X-1)\cong\F_q$ is involution
invariant, hence the equation reads $v^2\equiv 1$ modulo $(X-1)$,
namely $v\equiv\pm 1$ modulo $(X-1)$.
The family $\big((X-1)^i\big)_{0\le i\le n-1}$ is a basis of the
$\F_q$-vector space $\F_q[X]$, hence we write 
$v=\sum\limits_{k=0}^{n-1}v_k(X-1)^k\enspace,$
with $v_0=\pm 1$ and $v_k\in\F_q$. We compute the conjugates
$\wb{(X-1)^i}=(\wb{X}-1)^i$ of our
basis elements.
\begin{lemma}\label{conjugate}
For $0\le i\le n-1$, $(X-1)^i$ divides $(\wb{X}-1)^{i}$ and, more precisely:
$$(\wb{X}-1)^i=(-1)^i\sum_{k=0}^{n-i-1}\binom{n-i}{k}(X-1)^{k+i}\equiv
(-1)^{i}(X-1)^i\mod(X-1)^{i+1}\enspace.$$
\end{lemma}
\begin{proof}
Let $0\le i\le n-1$, then
$$(\wb{X}-1)^i=(X^{n-1}-1)^i=\big((1-X)X^{n-1}\big)^i=(-1)^i(X-1)^iX^{n-i}\enspace,$$
hence the equality, using Newton's formula for $X^{n-i}=(X-1+1)^{n-i}$.
\end{proof}
This result implies an important property for our approximation procedure.
\begin{lemma}
Let $1\le i\le \frac{n-1}{2}$, then
$$\left(\,v\wb v\equiv 1\mod(X-1)^{2i-1}\,\right)\Rightarrow\left(\,v\wb v\equiv 1\mod(X-1)^{2i}\,\right)\enspace.$$
\end{lemma}
\begin{proof}
Suppose the left hand side assertion is satisfied and write 
$$v\wb v\equiv 1+u(X-1)^{2i-1}\mod(X-1)^{2i}$$ 
for some $u\in\F_q$. Applying the involution shows that $(\wb{X}-1)^{2i}$ divides $v\wb
v-1-u(\wb{X}-1)^{2i-1}$, therefore
$$v\wb v\equiv 1+u(\wb{X}-1)^{2i-1}\mod(X-1)^{2i}\enspace,$$
thanks to Lemma \ref{conjugate}. We get:
$$0\equiv u\big((X-1)^{2i-1}-(\wb{X}-1)^{2i-1}\big)\equiv
2u(X-1)^{2i-1}\mod(X-1)^{2i}\enspace,$$
hence $u=0$.
\end{proof}
In particular we get that, if $v_0=\pm 1$, then $v\wb v\equiv
1\mod(X-1)^2$ for any value of $v_1\in\F_q$. We now need a formula for
the coefficients of $v$ of even positive index. 
\begin{lemma}\label{v2i}
Suppose $v\wb v\equiv 1\mod(X-1)^{2i}$ for some integer $1\le
i\le\frac{n-1}{2}$, then $v\wb v\equiv 1\mod(X-1)^{2i+1}$ if and only
if $v_{2i}$ satisfies (\ref{v}). 
\end{lemma}
\begin{proof}
Without any hypothesis on $v\wb v$, one checks using Lemma
\ref{conjugate} that:
$$v\wb v=\sum_{i=0}^{n-1}\left(\sum_{j=0}^{i}\sum_{k=0}^{j}(-1)^k\binom{n-k}{i-j}v_kv_{j-k}\right)(X-1)^i\enspace.$$
With our assumption on $v\wb v$, we get:
$$v\wb v\equiv 1+\sum_{j=0}^{2i}\sum_{k=0}^{j}(-1)^k\binom{n-k}{2i-j}v_kv_{j-k}\mod(X-1)^{2i+1}\enspace,$$
hence the result, noticing that $\binom{n}{2i}=0$ whereas $\binom{n}{0}=\binom{n-2i}{0}=1$.
\end{proof}
This ends the proof of Theorem \ref{thm:ramified}.
\end{proof}

%
\section{Experiments}\label{experiments}
\subsection{Algorithms}
Using MAGMA, we have implemented two algorithms based on the results of
this paper: the first finds a self-dual normal basis for a given
extension $\F_{q^n}/\F_q$ satisfying the existence conditions of
Theorem \ref{SDNBexistence} and such that the degree $n$ is either prime to the
characteristic or a power of it; the second computes the
orthogonal circulant group and uses it to construct all self-dual normal bases
of the extension from the former one, then selects those which have the
lowest complexity. Both these algorithms have a semi-simple and a
ramified version.

\subsubsection{Computation of a self-dual normal basis}
Our first algorithm permits us to find a self-dual normal basis for
somewhat large extensions. For example, one can find a self-dual
normal basis (of complexity $44\,431$) for $q=1009$ and $n=211$.
Here is the structure of this algorithm in the semi-simple case
$\gcd(n,q)=1$:  
\newcounter{savestep}
\begin{enumerate}[\sf Step 1.]
\item Compute the $q$-cyclotomic classes of the set $\{0,\ldots,n-1\}$.

\item Let $m$ be the size of the largest class (the class which
  contains~1) and choose $\zeta$ of order~$n$ in $\F_{q^m}$.

\item Build the matrices $F(\zeta)=(\zeta^{ij})_{1\leq i\leq j}$ and $F(\zeta^{-1})$.

\item Find a normal element~$\alpha$ in~$\F_{q^n}$. \textit{(This was already
  implemented in MAGMA, and uses methods which can be found in the
  book~\cite{Menezes})}. 

\item Compute $R\in\F_q[G]$ defined in Theorem \ref{vv}.  Using
  the matrix $F(\zeta)$, map $R$ to $R'=\FT(R)\in(\F_{q^m})^n$.

\item Use Lemmas~\ref{Erik3.5}, \ref{Erik3.6} and~\ref{Erik3.7} to
  find a solution $v'\in\mathop{\rm Im}\FT\subseteq(\F_q^m)^n$ of
  $v'\overline{v}'=R'$.  Bring back $v'$ to~$\F_q[G]$ using matrix
  $F(\zeta^{-1})$ to obtain $v$ such that $v\overline v=R$. Compute
  $w=v^{-1}$.

\item Compute and output $\gamma=w\circ\alpha$.
\setcounter{savestep}{\value{enumi}}
\end{enumerate}

In the odd characteristic, ramified case, we pick a normal
element~$\alpha$ in~$\F_{q^n}$ and compute $R\in\F_q[G]$; by
Proposition \ref{SDNBramified}, solving
the equation $v\overline v=R$ reduces to computing a square root of
$R$ in $\F_q[G]\simeq\F_q[X]/(X-1)^n$, which can be achieved by
computing a square root of $R$ modulo $X-1$ and then using Hensel
lifting.

\subsubsection{Computation of all self-dual normal bases of $\F_q^n$ over~$\F_q$}
The second algorithm can be used whenever the orthogonal circulant
group is not too large for an exhaustive enumeration, see Remark
\ref{size} and the tables in the next subsection. Here is its 
structure in the semi-simple case $\gcd(n,q)=1$: 
\begin{enumerate}[\sf Step 1.]
\setcounter{enumi}{\value{savestep}}
\item Use Proposition~\ref{ssc} to find generators (and their
  orders) of the group $U$ of solutions of $u\overline{u}=1$
  in~$\F_q[G]$ (this is actually done in the right hand side with
  elements of generators of $F_q^{m_k}$ where $m_k$ is the size of the
  cyclotomic class). 

\item For each $u$ in $U$ (elements of $U$ are enumerated using the
  generators found above), compute: the generator
  $\gamma=(uw)\circ\alpha$ of a self-dual normal basis, the
  multiplication-by-$\gamma$ matrix $\big(\tr(\gamma^{1+q^i+q^j})\big)_{i,j}$, and
  the complexity of~$\gamma$.  Update statistics accordingly (the best
  complexity found up to now, the list of best self-dual normal bases).

\item Finally, output the statistics (mainly the best complexity, and
  the number of times this complexity was achieved).
\end{enumerate}

In the ramified case, we list all the elements of $r\in\F_q[G]$
satisfying (\ref{r}), compute the associated $w$ as the proof 
of Theorem \ref{thm:MacWilliams} (\textit{i.e.} iteratively); the
group of solutions of $v\overline v=1$ consists of the elements $1+w$
obtained this way together with their opposites $-1-w$. We have each
of these elements act on the self-dual normal basis constructed above
and determine the complexity of the resulting self-dual normal basis.
\subsection{Tables}

The following tables show the complexity of the best self-dual normal
basis, obtained with the above algorithms, for some extensions.  We
give separate tables for extensions in characteristic~2 and for
extensions of small prime fields of odd characteristic. Blank entries
have not been computed since the cost of exhaustive enumeration grows
rapidly. 

\subsubsection{Even characteristic}\label{char2}

The lowest complexity for self-dual normal bases of extensions over
$\F_2$ is given in \cite[Table 5.1]{Jungnickel1993} for odd degree up
to $47$. With our method we were able to verify these values up to
$n=45$ (the computation for degree $45$ took approximately $25$ hours
on a $64$-bits Xeon quad core running at $2.33$ GHz).
We include our table for completeness.  
$$
  \vbox{\halign{%
      \strut\vrule\ $#$ \vrule width1.5pt&&\ $#$ \vrule\cr
      \noalign{\hrule}
            n      & 3 & 5 & 7 & 9 & 11 & 13 & 15 & 17 & 19 & 21 & 23 \cr
      \noalign{\hrule}
            \min   & 5 & 9 & 21 & 17 & 21 & 45 & 45 & 81 & 117 & 105 & 45 \cr
      \noalign{\hrule height1.5pt}
            n      & 25 & 27 & 29 & 31 & 33 & 35 & 37 & 39 & 41 & 43 & 45 \cr
      \noalign{\hrule}
            \min   & 93 & 141 & 57 & 237 & 65 & 69 & 141 & 77 & 81 & 165 & 153 \cr
      \noalign{\hrule}}}
$$
Note that \cite[Table 4]{Masuda_et_al} gives a minimal complexity of
$171$ for normal bases in degree $37$, where we find a self-dual normal
basis of complexity $141$, agreeing with Geiselmann
\cite[\textit{loc. cit.}]{Jungnickel1993}. Since only one digit
differs between these two results, we suspect that 
there may be 
a typo in \cite[\textit{loc. cit.}]{Masuda_et_al}. 

Using Lemma \ref{compositum}, one gets an upper bound for
the best self-dual normal complexities in even degree up to $n=90$,
using the fact that any element of $\F_4/\F_2$ of trace $1$ generates
an optimal self-dual normal basis (of complexity $3$). Comparing to
the results in \cite[Table 4]{Masuda_et_al} for $n$ up to $34$, we see
that this construction yields the best possible complexity in degrees
$10$, $22$ and $34$, and a reasonably good one in degrees $6$, $14$,
$18$, $26$ and $30$. 

We get optimal self-dual normal bases
in degrees $n=3$, $5$, $9$, $11$, $23$, $29$, $33$, $35$, $39$
and~$41$. We know by \cite[Corollary 3.6]{Mullin_et_al} that $2n+1$
has to be prime and $2$ of order $n$ or $2n$ modulo $2n+1$ for this to
happen, therefore we do not get optimal self-dual bases in degrees
$15$ and $21$, since $2$ is of order $5$ modulo $31$ and of order $14$
modulo $43$. 

We give also a table for other small even $q=2^r$.  Note that
$\alpha^{q^i}$ for $0\leq i\leq n-1$ generates the same normal basis
as $\alpha$, so the number of times the lowest complexity is obtained
is a multiple of~$n$.   When we found more than $n$ bases with the
lowest complexity, we indicate the multiplier between parentheses.
For example, we found 27 bases with complexity~45 for $q=8$ and~$n=9$.
$$
  \vbox{\halign{%
      \strut\vrule\ $#$ \vrule width1pt&&\ $#$ \vrule\cr
      \noalign{\hrule}
         q\backslash n & 3 & 5 & 7 & 9 & 11 & 13 & 15 & 17 & 19 & 21 & 23 & 25 \cr
      \noalign{\hrule height1pt}
            2    & 5 & 9 & 21 & 17 & 21 & 45 & 45 & 81 & 117(2) & 105 & 45 & 93 \cr
      \noalign{\hrule}
            4    & 5 & 9 & 21 & 17 & 21 & 45 & 45 & 81 & 117(2) & 105 & 45 & 93  \cr
      \noalign{\hrule}
            8    & 9(3) & 9 & 21 & 45(3) & 21 & 45 & 81(3) & 81 & & & & \cr
      \noalign{\hrule}
           16    & 5 & 9 & 21 & 17 & 21 & 45 & & & & & & \cr
      \noalign{\hrule}
           32    & 5 & 19(15) & 21 & 17 & 21 & & & & & & & \cr
      \noalign{\hrule}
           64    & 9(21) & 9  & 21 & 45(3) & & & & & & & & \cr
      \noalign{\hrule}
          128    & 5 & 9  & 37(98) & & & & & & & & & \cr
      \noalign{\hrule}
          256    & 5 & 9 & & & & & & & & & & \cr
      \noalign{\hrule}}}
$$

When
$\gcd(n,r)=1$ we always found the same best complexity for the 
extension $\F_{2^{rn}}$ over $\F_{2^r}$ as for the extension~$\F_{2^n}$
over~$\F_2$.
This observation is partially explained by
the following fact, which is also valid for odd $q$ (see \cite[Lemma
  4.2]{Menezes} for a partial proof). 
\begin{lemma}\label{EricLemma}
If $\alpha$ generates a self-dual normal basis of $\F_{q^n}$
over~$\F_q$, and $\gcd(n,r)=1$, then $\alpha$ generates a self-dual normal basis
of~$\F_{q^{rn}}$ over~$\F_{q^r}$, with the same complexity.
\end{lemma}

One easily checks that if an extension $\F_{q^n}/\F_q$ admits both
a self-dual normal basis and an optimal normal basis of type I (see
\cite{Gao+Lenstra}), then $q$ and $n$ have to be even, say $q=2^r$ and
$n=2m$, with $m$ odd and $2m+1$ prime. If this is the case, the
extension is the compositum of the fields $\F_{q^2}$ and $\F_{q^m}$,
each of which may admit an optimal self-dual normal basis or
not. Specifically, one can show that $\F_{q^2}/\F_{q}$ admits one
if and only if $r$ is odd, and that $\F_{q^m}/\F_{q}$ admits one if
$2$ is of order $m$ or $2m$ modulo $2m+1$ and $m$ is co-prime to
$r$. If all these conditions are satisfied, the self-dual normal basis
of $\F_{q^n}$ obtained by multiplying these two bases is, by Lemma 
\ref{compositum}, of complexity $3(2m-1)=3n-3$, which is also
the complexity of the dual basis of the optimal normal basis of
$\F_{q^n}$, see \cite[Theorem 5.4.10]{Jungnickel1993}
(\cite{Wan_Zhou} even shows that the dual of any basis which is
equivalent to the optimal one has complexity $3n-3$).  
This holds for instance for
the extensions of $\F_2$ of degrees $6$, $10$, $18$, $22$, $46$, ...,
those of $\F_8$ of degrees $10$, $22$, $46$, ... 

\subsubsection{Odd characteristic}

  Now we give the table showing some experiments for odd~$q$.  Here,
  the number of bases with least complexity is a multiple of~$2n$
  because $\pm\alpha^{q^i}$ for $0\leq i\leq n-1$ generates a normal
  basis with same complexity as the one generated by~$\alpha$.   The
  multipliers we indicate between parentheses, when we found more than
  $2n$ bases with lowest complexity, is relative to~$2n$. For example,
  we found $4\times2n=8n$ bases with complexity 51 for $q=13$
  and~$n=9$.  
$$
  \vbox{\halign{%
      \strut\vrule\ $#$ \vrule width1pt&&\ $#$ \vrule\cr
      \noalign{\hrule}
         q\backslash n & 3 & 5 & 7 & 9 & 11 & 13 & 15 & 17 & 19 & 21 & 23 & 25 \cr
      \noalign{\hrule height1pt}
           3  & \hbox{\bf7} & 13 & 25 & \hbox{\bf37} & 55 & 67 & -- & 91 & 172 & -- & 127 & 135 \cr
      \noalign{\hrule}
           5  &  6 & \hbox{\bf13} & 25 & 46 & 64(2) & 85 & -- & 157 & 153 & 150 & & \cr
      \noalign{\hrule}
           7  &  6 & 16 & \hbox{\bf19} & 41 & 61 & 96 & 87 & & & -- & & \cr
      \noalign{\hrule}
          11  &  6 & 13 & 25 & 52 & \hbox{\bf31} & 100 & 78 & & & & & \cr
      \noalign{\hrule}
          13  &  6 & 13 & 25 & 51(4) & 64 & \hbox{\bf37} & & & & & & \cr
      \noalign{\hrule}
          17  &  8 & 13 & 25 & 51(5) & 64 & 100 & & -- & & & & \cr
      \noalign{\hrule}
          19  &  8 & 13 & 31 & 51 & 67 & & & & -- & & & \cr
      \noalign{\hrule}}}
$$

Bold-face entries correspond to the best complexity in the case
when the degree $n$ is a power of the characteristic.  In this case
whenever $n$ is prime, the best complexity is $3n-2$, and is obtained
with the basis exhibited in~\cite[Theorem 5.3]{Blake_Gao_Mullin}. This
basis is rather explicit since generated by the root of a
trinomial, yielding a very interesting family of self-dual normal
bases, of complexity fairly close to the optimal one.

We have made no computation for ``mixed degree'' $n=n_1p^e$ with
$\gcd(n_1,p)=1$, $n_1>1$ and $e>0$, but one gets an upper
bound for the lowest complexity in that case by multiplying the lowest
complexity in degree $n_1$ by that in degree $p^e$, thanks to Lemma
\ref{compositum}. For instance, the best complexity for $q=5$ and
$n=15$ is at most $6\cdot 13=78$. Note 
that when $n=\ell\ell'$ for prime numbers $\ell\not=\ell'$, both
different from $p$, the best complexity for the compositum is not
necessarily the product of those for degrees $\ell$ and $\ell'$
extensions ($n=15$, $q=7$); however it can be so ($n=15$, $q=11$; $n=21$,
$q=5$).

In the semi-simple case, we also computed the best complexity for some
odd non prime values $q=p^r$, which do not appear in this table.  
When $\gcd(n,r)=1$ we always found the same best complexity for the 
extension $\F_{q^n}$ over $\F_q$ as for the extension~$\F_{p^n}$
over~$\F_p$, as well as the same multiplier for the number of bases
with the best complexity (as in the even characteristic case).



In odd characteristic, the only exhaustive search for lowest
complexities among normal bases we are aware of is in
\cite{Blake_Gao_Mullin}, over prime base fields. 
The lowest complexity for self-dual normal bases is the same as the
one they obtain for normal bases
when $n=3$ and $q=7$ or $13$; slightly larger when $n=3$ and
$q=19$ ($8$ instead of $6$) and when $n=5$ and $q=11$ ($13$ instead of
$12$). Note that in this last case, Liao and Feng give in 
\cite[Example 2]{Liao_Feng} a construction of a normal basis with
minimal complexity $12$, using Gauss periods, whose dual basis has
complexity $13$. Their construction remains valid when replacing the
base field $\F_{11}$ by an extension of degree prime to $5$. 

\subsection{Conclusion}
Our algorithms enable us to compute the minimal complexity for
self-dual normal bases in various extensions of finite
fields, including some for which the
exhaustive enumeration of normal bases would not be reasonable.
In odd characteristic, the lowest complexities we obtain are
either the same as or close to that obtained in
former computations on normal bases using theoretical constructions or
exhaustive search, analogously to what could already be observed in even
characteristic. However the cost of 
the exhaustive search of all self-dual normal bases (once one has been
constructed) is still a limitation of this method. In order to make
self-dual normal bases practical,
it would thus be desirable to find
a direct construction of those with low complexity.

A striking fact when looking at the tables above is the
repetition of values along columns, albeit with some exceptions. We
have a partial explanation for this phenomenon, that may also help in
achieving the former goal, in terms of global considerations of
cyclotomic extensions of the rationals generated by $n^2$-th roots of
unity, where $n$ is a prime. A known construction yields a global
self-dual normal basis generator $\alpha_n$ such that, for any prime
$p\not=n$ which does not split in the considered extension, the
residue modulo $p$ of $\alpha_n$ is a candidate for a best complexity
basis for $\F_{p^n}/\F_p$. We hope to give full details about this 
construction in a future paper.

\begin{aknowledgments}
The authors would like to thank the two anonymous referees for their
valuable and insightful remarks and advice.
\end{aknowledgments}

\bibliography{bib}

\begin{thebibliography}{10}
\expandafter\ifx\csname url\endcsname\relax
  \def\url#1{\texttt{#1}}\fi
\expandafter\ifx\csname urlprefix\endcsname\relax\def\urlprefix{URL }\fi

\bibitem{Ash_Blake_Vanstone}
D.~W. Ash, I.~F. Blake, S.~A. Vanstone, Low complexity normal bases, Discrete
  Appl. Math. 25~(3) (1989) 191--210.

\bibitem{BG}
T.~Beth, W.~Geiselmann, Selbstduale {N}ormalbasen \"uber {${\rm GF}(q)$}, Arch.
  Math. (Basel) 55~(1) (1990) 44--48.

\bibitem{Blake_Gao_Mullin}
I.~F. Blake, S.~Gao, R.~C. Mullin, Normal and self-dual normal bases from
  factorization of {$cx^{q+1}+dx^q-ax-b$}, SIAM J. Discrete Math. 7~(3) (1994)
  499--512.

\bibitem{ByrdVaughan}
K.~A. Byrd, T.~P. Vaughan, Counting and constructing orthogonal circulants, J.
  Combinatorial Theory Ser. A 24~(1) (1978) 34--49.

\bibitem{Christopoulou_et_al}
M.~Christopoulou, T.~Garefalakis, D.~Panario, D.~Thomson, The trace of an
  optimal normal element and low complexity normal bases, Des. Codes Cryptogr.
  49~(1-3) (2008) 199--215.

\bibitem{Couveignes_Lercier}
J.-M. Couveignes, R.~Lercier, Elliptic periods for finite fields, Finite Fields
  Appl. 15~(1) (2009) 1--22.

\bibitem{Gao}
S.~Gao, Normal {B}ases {O}ver {F}inite {F}ields, {PhD} in {C}ombinatorics and
  {O}ptimisation, University of {W}aterloo, Waterloo, {O}ntario, {C}anada
  (1993).

\bibitem{Gao+Lenstra}
S.~Gao, H.~W. Lenstra, Jr., Optimal normal bases, Des. Codes Cryptogr. 2~(4)
  (1992) 315--323.

\bibitem{Gao_et_al}
S.~Gao, J.~Von Zur~Gathen, D.~Panario, V.~Shoup, Algorithms for exponentiation
  in finite fields, J. Symbolic Comput. 29~(6) (2000) 879--889.

\bibitem{Jungnickel1993}
D.~Jungnickel, Finite fields, Bibliographisches Institut, Mannheim, 1993,
  structure and arithmetics.

\bibitem{J-B-G}
D.~Jungnickel, T.~Beth, W.~Geiselmann, A note on orthogonal circulant matrices
  over finite fields, Arch. Math. (Basel) 62~(2) (1994) 126--133.

\bibitem{JMV}
D.~Jungnickel, A.~J. Menezes, S.~A. Vanstone, On the number of self-dual bases
  of {${\rm GF}(q^m)$} over {${\rm GF}(q)$}, Proc. Amer. Math. Soc. 109~(1)
  (1990) 23--29.

\bibitem{LempelWeinberger}
A.~Lempel, M.~J. Weinberger, Self-complementary normal bases in finite fields,
  SIAM J. Discrete Math. 1~(2) (1988) 193--198.

\bibitem{Liao_Feng}
Q.~Liao, K.~Feng, On the complexity of the normal bases via prime {G}auss
  period over finite fields, J. Syst. Sci. Complex. 22~(3) (2009) 395--406.

\bibitem{Liao_Sun}
Q.~Y. Liao, Q.~Sun, Normal bases and their dual-bases over finite fields, Acta
  Math. Sin. (Engl. Ser.) 22~(3) (2006) 845--848.

\bibitem{MacWilliams}
F.~J. MacWilliams, Orthogonal circulant matrices over finite fields, and how to
  find them., J. Combinatorial Theory Ser. A 10 (1971) 1--17.

\bibitem{Masuda_et_al}
A.~M. Masuda, L.~Moura, D.~Panario, D.~Thomson, Low complexity normal elements
  over finite fields of characteristic two, IEEE Trans. Comput. 57~(7) (2008)
  990--1001.

\bibitem{Menezes}
A.~J. Menezes, I.~F. Blake, S.~Gao, R.~C. Mullin, S.~A. Vanstone, T.~Yaghoobian
  (eds.), Applications of finite fields., Kluwer Academic Publishers, 1993.

\bibitem{Mullin_et_al}
R.~C. Mullin, I.~M. Onyszchuk, S.~A. Vanstone, R.~M. Wilson, Optimal normal
  bases in {${\rm GF}(p^n)$}, Discrete Appl. Math. 22~(2) (1988/89) 149--161.

\bibitem{Nogami_et_al}
Y.~Nogami, H.~Nasu, Y.~Morikawa, S.~Uehara, A {M}ethod for {C}onstructing a
  {S}elf-{D}ual {N}ormal {B}asis in {O}dd {C}haracteristic {E}xtension
  {F}ields, Finite Fields Appl. 14 (2008) 867--876.

\bibitem{Pickett2}
E.~J. Pickett, Construction of {S}elf-{D}ual {I}ntegral {N}ormal {B}ases in
  {A}belian {E}xtensions of {F}inite and {L}ocal {F}ields, Int. J. Number
  Theory 6~(7) (2010) 1565--1588.

\bibitem{Poli}
A.~Poli, Constructing {SCN} bases in characteristic {$2$}, IEEE Trans. Inform.
  Theory 41~(3) (1995) 790--794.

\bibitem{Wan_Zhou}
Z.-X. Wan, K.~Zhou, On the complexity of the dual basis of a type {I} optimal
  normal basis, Finite Fields Appl. 13~(2) (2007) 411--417.

\bibitem{Wang}
C.~C. Wang, An {A}lgorithm to {D}esign {F}inite {F}ield {M}ultipliers {U}sing a
  {S}elf-{D}ual {N}ormal {B}asis, IEEE Transactions on Computers 38~(10) (1989)
  1457--1460.

\end{thebibliography}

\end{document}